\documentstyle[12pt]{article}
\newcommand{\braid}[3]{{#1}$\lower4pt\hbox{$\oo\atop\raise4pt
           \hbox{$\scriptscriptstyle {#3} $}$}${#2}}
\newcommand{\twist}[3]{{#1}${\,\scriptscriptstyle {#3}}\atop\raise9pt
           \hbox{$\scriptstyle\oo$} ${#2}}

\newcommand{\be}{\begin{eqnarray}}
\newcommand{\ee}{\end{eqnarray}}
\newcommand{\n}{\nonumber }
\newcommand{\oo}{\otimes}

\newcommand{\la}{\lambda}

\def\C{{\mathchoice {\setbox0=\hbox{$\displaystyle\rm C$}\hbox{\hbox
to0pt{\kern0.4\wd0\vrule height0.9\ht0\hss}\box0}}
{\setbox0=\hbox{$\textstyle\rm C$}\hbox{\hbox
to0pt{\kern0.4\wd0\vrule height0.9\ht0\hss}\box0}}
{\setbox0=\hbox{$\scriptstyle\rm C$}\hbox{\hbox
to0pt{\kern0.4\wd0\vrule height0.9\ht0\hss}\box0}}
{\setbox0=\hbox{$\scriptscriptstyle\rm C$}\hbox{\hbox
to0pt{\kern0.4\wd0\vrule height0.9\ht0\hss}\box0}}}}

\newtheorem{theorem}{Theorem}
\begin{document}
\begin{titlepage}

\begin{center}
{\Large \bf On twisting solutions to the Yang-Baxter equation
 }
\end{center}
\vspace{0.5cm}
\begin{center}
 {\bf P. P. Kulish\footnote{Partially supported by the RFFI
 grant 98-01-00310.}}
\end{center}
\vspace{0.5cm}
\begin{center}
St.Petersburg Department of the Steklov
Mathematical Institute,\\
Fontanka 27, St.Petersburg, 191011,
Russia \\ (kulish@pdmi.ras.ru)
\end{center}
\vspace{0.5cm}
\begin{center}
 {\bf A. I. Mudrov}
\end{center}
\vspace{0.5cm}
\begin{center}
Department of Theoretical Physics,
Institute of Physics, St.Petersburg State University, Ulyanovskaya 1,
St.Petergof, St.Petersburg, 198904, Russia\\
(aimudrov@dg2062.spb.edu)
\end{center}
\vspace{1cm}
\begin{center}
{\bf Abstract}
\end{center}
{\footnotesize
Sufficient conditions for an invertible two-tensor
$F$ to relate two solutions to the Yang-Baxter equation
via the transformation $R\to F^{-1}_{21} R F$
are formulated. Those conditions
include relations arising from twisting of certain
quasitriangular bialgebras.
}
\vspace{0.5in}
\end{titlepage}
\section{Introduction}
The twist procedure for (quasi)-Hopf algebras developed by
Drinfeld  \cite{D1,D2} (see also \cite{R,RSTS})
allows to deform the coproduct via a similarity transformation with the
multiplication unchanged. Twist finds its applications in solvable
models and noncommutative geometry because it appears to be
very friendly to all the algebraic properties of a given Hopf
algebra, including the quasitriangular structure and the structure of
modules.
However, to find the explicit form of an element
${\cal F}$ realizing interrelations between twisted and untwisted objects
is as difficult problem as that of evaluating universal
R-matrices ${\cal R}$. On the other hand, most applications of quantum
groups  employ their particular matrix representations,
and in practice one deals with matrix solutions to the Yang-Baxter
equation (YBE) rather than the universal ones. The FRT algorithm \cite{RTF}
yields the recipe how, starting from finite-dimensional solutions to
YBE, to build both quantum groups and universal R-matrices \cite{D3}.
The latter is possible due to the factorization
property of ${\cal R}$  with respect to the coproduct and the famous
fusion procedure
\cite{KRS}. It is worth to note that the factorization property
is virtually the only tool known for building universal twisting
elements as well as their matrix realizations
\cite{R,RSTS},\cite{K1}--\cite{KLM}. The finite-dimensional
(matrix) version  of twisting procedure was considered in the general
setting in \cite{M}, but up to now there are no general criteria, except
already mentioned factorization properties, for a matrix two-tensor $F$
to define twist of
a given Hopf algebra. On the other hand, the problem of transforming
solutions to the Yang-Baxter equation $R\to F^{-1}_{21} R F$
makes sense by itself, regardless of the possibility of
expanding $F$ to the universal element ${\cal F}$. In the present paper we
formulate sufficient conditions guaranteeing that the transformation
of concern should provide a new solution to YBE. Quite amazingly,
they involve an invertible three-tensor which itself drops from
YBE but ensures its fulfillment.

  The article is organized as follows. Section II is devoted to
the transformation of solutions to YBE which we, by the analogy with
that arising within Drinfeld's theory, call twist.
The relation between twist of bialgebras and twist of R-matrices
is discussed in Section III. Section IV demonstrates some examples
when twist of a matrix solution to YBE can be extended to the global
twist of the quantum algebra dual to the corresponding FRT quantum
semi-group. In Conclusion we discuss possible applications of the
results obtained.

\section{Twists of R-matrices}
Throughout the paper ${\bf R}$ will denote an associative algebra
with unit over a field $k$.
The main result of the present communication is given by the following
assertion.
\begin{theorem}
\label{Th2}
Let  $R\in {\bf R}^{\oo 2}$ be a
solution to the Yang-Baxter equation
$$
R_{12}R_{13}R_{23}=R_{23}R_{13}R_{12},
$$
and invertible elements $F \in {\bf R}^{\oo 2} $,
$\Phi \in {\bf R}^{\oo 3} $,
and $\Psi \in {\bf R}^{\oo 3} $ fulfill the following conditions
\be
\Phi_{123}F_{12}=\Psi_{123}F_{23}
\label{cond1}
\ee
\be
R_{12}\Phi_{123}=\Phi_{213}R_{12}
\label{cond2}
\ee
\be
R_{23}\Psi_{123}=\Psi_{132}R_{23}
\label{cond3}
\ee
Then
\be
\tilde R =\tau(F^{-1})R F,
\label{twist}
\ee
is a solution to the Yang-Baxter equation too ($\tau$ is the
permutation of the tensor components).
\end{theorem}
Transformation (\ref{twist}) is called  {\em twist} of a solution to YBE.
For a matrix ring ${\bf R}$ this turns into a similarity transformation of
the braid matrix $\hat R=PR$, with $P$ representing the permutation
operation in the corresponding vector space.
In fact, twist is determined by a pair $(F,G)$, $G\in {\bf R}^{\oo 3}$, such
that $\Phi=G\bar F_{12}$ and $\Psi=G\bar F_{23}$
obey the conditions of the theorem
(to make formulas more readable
we denote the inverse by bar).

First let us prove the equality
\be
\bar F_{12} \bar \Psi_{312}R_{13} R_{23}
 \Phi_{123}F_{12} = \tilde R_{13} \tilde R_{23}.
\label{ax}
\ee
Indeed, using conditions (\ref{cond1}--\ref{cond3}) along with
the definition (\ref{twist}) we find
\be
\bar F_{12} \bar \Psi_{312}R_{13} R_{23}\Phi_{123}F_{12}
&=&
\bar F_{12} \bar \Psi_{312}R_{13} R_{23}\Psi_{123}F_{23}
\qquad
\qquad
\qquad
\qquad
\qquad
\qquad
\n\\&=&
\bar F_{12} \bar \Psi_{312}R_{13} \Psi_{132} R_{23} F_{23}
\n\\&=&
\bar F_{12} \bar \Psi_{312}R_{13} \Psi_{132} F_{32} \tilde R_{23}
\n\\&=&
\bar F_{12} \bar \Psi_{312}R_{13} \Phi_{132} F_{13} \tilde R_{23}
\n\\&=&
\bar F_{12} \bar \Psi_{312} \Phi_{312} R_{13}  F_{13} \tilde R_{23}
\n\\&=&
\bar F_{12} \bar \Psi_{312} \Phi_{312} F_{31} \tilde R_{13}\tilde R_{23}
\n\\&=&
\bar F_{12} \bar \Psi_{312} \Psi_{312} F_{12} \tilde R_{13}\tilde R_{23}
\n\\&=&
 \tilde R_{13} \tilde R_{23}.
\n
\ee
Now, taking into account the auxiliary identity (\ref{ax}), we have
\be
\tilde R_{12}\tilde R_{13} \tilde R_{23} {\tilde R_{12}}^{-1}
&=&
 \bar F_{21}R_{12} F_{12}
 \bar F_{12} \bar \Psi_{312}R_{13} R_{23}\Phi_{123}F_{12}
 \bar F_{12}\bar R_{12} F_{21}
\n\\&=&
 \bar F_{21}R_{12}
 \bar \Psi_{312}R_{13} R_{23}\Phi_{123}
 \bar R_{12} F_{21}
\n\\&=&
 \bar F_{21} \bar \Psi_{321}R_{12}R_{13} R_{23}\bar R_{12}\Phi_{213}  F_{21}
\n\\&=&
 \bar F_{21} \bar \Psi_{321} R_{23} R_{13}\Phi_{213}  F_{21}
\n\\&=&
 \tau_{12}
 (\bar F_{12} \bar \Psi_{312}R_{13} R_{23}\Phi_{123}F_{12})
\n\\&=&
 \tau_{12}
 (\tilde R_{13}\tilde R_{23})= \tilde R_{23}\tilde R_{13},
\ee
as required.

Theorem \ref{Th2} can be understood within the framework of the
bialgebra twist theory in its matrix formulation rendered in some detail
in the next section. Although twist of an R-matrix might not be extended to
the global bialgebra twist as discussed later on,  it possesses many
familiar features, for example, the composition property.
For a pair $(F,G)$ the inverse is $(\bar F,\bar G )$.
For two pairs $(F,G)$ and $(F',G')$, where
the latter is defined through the new R-matrix $\tau(\bar F)R F$, there
exists their composition $(F F', GG')$.
And, finally, $(e^2,e^3)$ (the units in the tensor
square and cube of ${\bf R}$) realizes the identical transformation.
Another similarity with the bialgebra twist is the existence of
the set of gauge transformations of the twisting pair as
will be shown in the next section.

\section{On the global twist and braid group representation equivalence}
Although the most evident applications of the observation made
in the previous section can be relevant to finite-dimensional matrix
rings ${\bf R}$, usually those of fundamental representations
of Hopf algebras of interest, Theorem \ref{Th2} holds for any ${\bf R}$.
So, one can consider ${\bf R}= $Mat$(N)[[\la,\la^{-1}]]$ and R-matrices
depending on the spectral parameter $\la$. As another example,  let us take a
quasitriangular bialgebra ${\cal H}$ as ${\bf R}$
and an element ${\cal F}\in{\cal H}^{\oo2}$ satisfying the twist
equation \cite{D2}
\be
 (\Delta\oo id)({\cal F}){\cal F}_{12}= (id\oo\Delta)({\cal F}){\cal F}_{23}.
\label{TE}
\ee
Then, one can put $F={\cal F}$ and $G$ to be the expression on either side of
(\ref{TE}).

To explain the result obtained we shall use the formalism dual to
the FRT algorithm of constructing quantum semi-groups. Recall
that the tensor bialgebra $T({\bf R})$ over ${\bf R}$ is introduced as
the direct sum of ideals $T({\bf R})=\sum^{\infty}_{n=0} {\bf R}^{\oo n}$,
where  ${\bf R}^0$ is isomorphic to $k$,  the field of scalars.
The unit in $T({\bf R})$ is represented by the sum of idempotents
$1=\sum_{n\geq 0} e^n$, the units in ${\bf R}^{\oo n}$,
respectively. Multiplication by $e^n$ realizes the
projection homomorphism $T({\bf R})\to {\bf R}^{\oo n}$, and
for $n=0$ this coincides with the counit mapping to $k$, the coproduct
being introduced on the basis elements $x^{i_1..i_n}\in {\bf R}^{\oo n}$ as
\be
  \Delta (x^{i_1..i_n}) &=& e^0\oo x^{i_1..i_n}+...+x^{i_1..i_k}
  \oo x^{i_{k+1}..i_n}+...+x^{i_1..i_n}\oo e^0
  .\n
\ee
The principal feature of $T({\bf R})$
is that for any bialgebra ${\cal H}$ and a representation
$\rho\colon {\cal H} \to {\bf R}$ there is the unique extension
to the homomorphism of bialgebras ${\cal H} \to T({\bf R})$.
It is built by means of the multiple coproduct
$\Delta^n\colon {\cal H}\to{\cal H}^{\oo n}$ defined for
$n=0$ as the counit, for $n=1$ as the identical mapping, and
for higher $n$ it is $\Delta^2\equiv \Delta,\>$
$\Delta^3\equiv (\Delta\oo id)\circ\Delta$, and so on.
Then the homomorphism of concern is specified by the
mappings $(\rho^{\oo n}\circ\Delta^n )\colon {\cal H}\to {\bf R}^{\oo
n}$.

For a given solution to the Yang-Baxter equation $R \in {\bf R}^{\oo 2}$,
one defines a subalgebra ${\cal U}=\sum^{\infty}_{n=0} {\cal U}^n$, where
${\cal U}^0= k$, ${\cal U}^1={\bf R}$, and
${\cal U}^n=\{z|z\in {\bf R}^{\oo n}, R_{ii+1} z =
\tau_{ii+1}(z) R_{ii+1},0<i< n\}$
($\tau_{ii+1}$ is the permutation between $i$-th and $i+1$-th sites).
Such tensors are called $R$-symmetric, and in the case of matrix rings
they just commute with the braid matrix $\hat R$.
\begin{theorem}
${\cal U}$ is a quasitriangular sub-bialgebra in $T({\bf R})$.
\end{theorem}
It follows immediately from the definition that
${\cal U}$ is a sub-bialgebra indeed . Its universal
R-matrix is decomposed into the sum of its
${\bf R}^{\oo m}\oo {\bf R}^{\oo n}$-components $R^{m,n}$;
for $ m n =0$ it is just $e^m\oo e^{n}$, the unit of
${\bf R}^{\oo m}\oo {\bf R}^{\oo n}$, and if $n=n'$ and $m$ are both
non-zero, one has
\be R^{m,n'}=(R_{1n'}\ldots R_{11'})(R_{2n'}\ldots
 R_{21'})\ldots (R_{mn'}\ldots R_{m1'}), \label{RM} \ee
where primes
mark indices numbering ${\bf R}$-factors in the second tensor
component. Note that the bialgebra ${\cal U}$ is dual to the
quantum semi-group ${\cal A}_R$ generated by basis elements of the
linear space ${\bf R}^*$ with the RTT relations imposed.

Usually, one is interested in Hopf algebras, which
require additional relations  of the quantum determinant type
imposed on the generators of the quantum semi-group. Such
relations eliminate just few degrees of freedom while dramatically
complicate algebraic structure mixing the homogeneous components.
So, we prefer to work with bialgebras, in the dual sector
represented by the direct sum of their ideals ${\cal U}=\sum_{n\geq
0}{\cal U}\cap {\bf R}^{\oo n}$.

The component representation of the twist equation (\ref{TE})
in $T({\bf R})$ reads
\be
F^{m+n,k}(F^{m,n}\oo e^k)= F^{m,n+k}(e^m\oo F^{n,k}),
\label{TE1}
\ee
with $F^{i,j}$ being the images of the universal twisting element
${\cal F}$ in
${\bf R}^{\oo i}\oo {\bf R}^{\oo j}\subset T({\bf R})^{\oo 2}$.
It is equal to
$F^{m,n}=(\rho^{\oo m}\circ\Delta^m\oo\rho^{\oo n}\circ\Delta^n)({\cal F})$,
where $\rho$ is the representation of ${\cal H}$ in ${\bf R}$.
Given a solution to (\ref{TE}),
for any quasitriangular bialgebra one can construct
twisted quasitriangular bialgebra with the new universal R-matrix
${\cal R}=\tau({\cal F})^{-1}{\cal R}{\cal F}.$  An interesting implication
of the equivalent system (\ref{TE1}) is that there are no closed conditions
on $F^{1,1}$ directly involved in deformation of
$R=R^{1,1}=(\rho\oo\rho)({\cal R})$, which is a matrix solution to
YBE.  Actually, to obtain new, twisted solutions to the Yang-Baxter
equation there is no need to satisfy the whole set of equations
(\ref{TE1}) recovering  the universal element ${\cal F}$, it
is sufficient to restrict the study by only the small part of them.
This is the observation which underlines Theorem~ \ref{Th2}.

Let us investigate the question when a twist of an R-matrix can be extended
to the twist of the entire bialgebra ${\cal U}$. Having introduced tensors
$\Omega^2=F$, $\Omega^3=G$, one can see that they satisfy
the equalities
\be
R_{ii+1} \Omega^n = \tau_{ii+1}(\Omega^n) \tilde R_{ii+1},\quad 0<i<n
\label{braid}
\ee
for $n=2,3$.
 If ${\bf R}$ is a matrix ring, this establishes
a local isomorphism of the braid group $B_3$ local representations
specified by the matrices $R$ and $\tilde R$.
\begin{theorem}
\label{Th3}
The pair $(\Omega^2,\Omega^3)$ is extended to the twist
of the bialgebra ${\cal U}$ if and only if for each  $n>3$ there exists
an invertible element $\Omega^n\in {\bf R}^{\oo n}$ fulfilling
(\ref{braid}). Twisting element is uniquely defined up to
an isomorphism via the formula
\be
F^{m,n}=\Omega^{m+n}(\bar\Omega^{m}\oo\bar\Omega^{n}),\quad m,n\geq 0,
\label{TE2}
\ee
where for $i=0,1$ we set $\Omega^i=e^i$, the units of
${\bf R}^{\oo i}$.

\end{theorem}

It is easy to see that $F^{m,n}$ introduced according to (\ref{TE2})
satisfy (\ref{TE1}) and indeed lie in ${\cal U}^{\oo 2}$ (their each
component
is $R$-symmetric). Let us prove the converse. Given the universal twisting
element, define $\Omega^n$ for $n>2$ as the product
\be
\Omega^n=
F^{1,n-1}(e^1\oo F^{1,n-2})...(e^{n-2}\oo F^{1,1})
=F^{1,n-1}(\Omega^1\oo \Omega^{n-1})
.
\label{Int}
\ee
We are going to state (\ref{TE2}) and that would
evidently be enough because then we can employ the induction method
and the $R$-symmetry of the elements $F^{m,n}$.
Conditions (\ref{TE1}) hold if one of the numbers $m$ and $n$ are zero. They
are also true by construction for $m=1$ and any $n$. Then, for $m\geq 1$ one
has
\be
\Omega^{1+m+n}&=&
F^{1,m+n}(\Omega^{1}\oo \Omega^{m+n})=
F^{1,m+n} \biggl(e^{1}\oo\Bigl(F^{m,n}(\Omega^{m}\oo\Omega^{n})\Bigr)
\biggr)=
\n\\&=&
F^{1,m+n} (e^{1}\oo F^{m,n})(e^{1}\oo\Omega^{m}\oo\Omega^{n})=
\n\\&=&
F^{1+m,n} (F^{1,m}\oo e^{n})(e^{1}\oo\Omega^{m}\oo\Omega^{n})=
\n\\&=&
F^{1+m,n}\Bigl(F^{1,m} (e^{1}\oo\Omega^{m})\oo\Omega^{n}\Bigr)=
F^{1+m,n}(\Omega^{1+m}\oo\Omega^{n})
\n
\ee
by induction. Thus, as subalgebras in $T({\bf R})$, ${\cal U}$ and its
twisted counterpart $\tilde {\cal U}$
are related by the similarity transformation with
the element $\Omega=\sum_{n\geq 0}\Omega^n$ and {\em vice versa}.
This immediately implies the uniqueness of the global twist
because two different $\Omega$'s are linked via an $R$-symmetric element
$u$ which, by definition, belongs to ${\cal U}$ (that is also a
manifestation
of the twist composition property). It realizes the inner automorphism
$h\to u^{-1} h u$ leading to the transformation
${\cal F}\to \Delta(u^{-1}){\cal F}(u\oo u)$.
If it happens so that given $\Omega^i$, $i=2,3$,
cannot be expanded to a universal element ${\cal F}$,
yet there are
gauge transformations of the element $F$ leading to trivial
or isomorphic deformations of the R-matrix.  For every invertible
$R$-symmetric $u^i\in {\bf R}^{\oo i}$, $i=1,2,3$,
substitution
$(\Omega^2,\Omega^3) \to \Bigl( u^2 \Omega^2 (u^1\oo u^1),
u^3 \Omega^3 (u^1\oo u^1\oo u^1)\Bigr)$,
results in the similarity transformation
$\tilde R\to (\bar u^1\oo \bar u^1) \tilde R (u^1\oo u^1)$.

We conclude this section with the remark that in the case of quasitriangular
bialgebra ${\cal H}$ admitting twist with the element ${\cal F}$ there is
the abstract form of $F^{m,n}$ belonging to ${\cal H}^{\oo
m}\oo{\cal H}^{\oo n}$.  It is built with the help of the multiple coproduct
applied to the components of the twisting element:
${\cal F}^{m,n}=(\Delta^m\oo\Delta^n)({\cal F})$.  Formula (\ref{TE2})
then gives the abstract element $\Omega$ intertwining ${\cal R}$- and
$\tilde {\cal R}$-symmetric tensors in ${\cal H}^{\oo n}$ and in
algebraically isomorphic $\tilde {\cal H}^{\oo n}$. Element $\Omega$
appeared in \cite{MS} as the necessary condition for the global twist
factorization of the unitary universal R-matrix and was applied to the
$XXZ$-model of spin $\frac{1}{2}$ chain described by the trigonometric
solution to the Yang-Baxter equation. The relation to the symmetric
group $S_n$ representations was discussed there as well.

\section{Factorization of twisting elements}
In the present section we give some illustrations
to the constructions considered above.
Given a bialgebra ${\cal H}$ and a solution ${\cal F}\in {\cal H}\oo{\cal H}$
to the pair of equations
\be
(id \oo\Delta)({\cal F})={\cal F}_{13}{\cal F}_{12},
\label{FC11}
\ee
\be
(\Delta\oo id)({\cal F})={\cal F}_{13}{\cal F}_{23},
\label{FC12}
\ee
satisfying the additional condition
\be
{\cal F}_{12} {\cal F}_{23}={\cal F}_{23}{\cal F}_{12},
\label{FC13}
\ee
it is possible to twist ${\cal H}$ by ${\cal F}$. In a matrix
representation, the element ${\cal F}$ is given by its tensor components
decomposed into the products
$$ F^{m,n'}=(F_{1n'}\ldots
F_{11'})(F_{2n'}\ldots F_{21'})\ldots (F_{mn'}\ldots F_{m1'})$$
(cf. the notation in formula (\ref{RM})).
This expression is exactly the
same as for the universal R-matrix (\ref{RM}), and that is due to
the factorization conditions (\ref{FC11}),  (\ref{FC12}) similar
to those held for universal R-matrices.  The homogeneous components of
the global intertwiner $\Omega$ are evaluated using (\ref{Int}) and
(\ref{FC11}):
$$ \Omega^n = (F_{1\ n}\ldots F_{1 2})( F_{2\ n}\ldots F_{2 3} )
   \ldots (F_{n-1\ n}).  $$
The most natural situation for such
twists appears when the bialgebra ${\cal H}$ is isomorphic to the
tensor product of its sub-bialgebras ${\cal A}$ and ${\cal B}$ and
${\cal F}$ actually belongs to ${\cal A}\oo{\cal B}$ \cite{RSTS}.
In the matrix language (\ref{FC11}-\ref{FC13}) read
\be
R_{23}F_{13}F_{12}=F_{12}F_{13}R_{23},
\label{Split11}
\ee
\be
R_{12}F_{23}F_{13}=F_{13}F_{23}R_{12}.
\label{Split12}
\ee
\be
F_{12}F_{23}=F_{23}F_{12}.
\label{Split13}
\ee
Conversely, each solution to the system (\ref{Split11}--\ref{Split13})
generates a twist of the R-matrix extended to the global twist of
${\cal U}$. Twisting elements fulfilling (\ref{Split11}--\ref{Split13})
were used for explaining Fronsdal-Galindo deformation of the standard
Drinfeld-Jimbo quantum groups ${\cal U}_q(sl(2N+1))$ \cite{JC}.

Another possible factorizations of the twisting element with respect to
the coproducts (in this case, twisted and non-twisted ones) are \cite{M,KLM}
\be
(id \oo\Delta)({\cal F})={\cal F}_{12}{\cal F}_{13}.
\label{FC21}
\ee
\be
(\tilde\Delta\oo id)({\cal F})={\cal F}_{13}{\cal F}_{23}.
\label{FC22}
\ee
The system of equations (\ref{FC21}), (\ref{FC22}),
and  (\ref{TE}) is determined by its any pair. Conditions
(\ref{FC21}--\ref{FC22}) are the generalization of
Reshetikhin's twist, in which (\ref{FC22}) is substituted by
$(\Delta\oo id)({\cal F})={\cal F}_{23}{\cal F}_{13}$ and the
Yang-Baxter relation ${\cal F}_{12}{\cal F}_{13}{\cal F}_{23}=
{\cal F}_{23}{\cal F}_{13}{\cal F}_{12}$ (originally there were
some excessive additional conditions which were loosened later in \cite{K1}).
Matrix version of (\ref{FC21}--\ref{FC22}) reduces to
\be
R_{23}F_{12}F_{13}=F_{13}F_{12}R_{23},
\label{Split21}
\ee
\be
\tilde{R}_{12}F_{13}F_{23}=F_{23}F_{13}\tilde{R}_{12}.
\label{Split22}
\ee
Indeed, as was shown in \cite{M}, any $F$ fulfilling
(\ref{Split21}--\ref{Split22}) defines the global twist
possessing (\ref{FC21}--\ref{FC22}). Again, using factorization
(\ref{FC21}) we find
$$
\Omega^n = (F_{1\ 2}\ldots F_{1\ n}) (F_{2\ 3}\ldots F_{2\ n})
\ldots  (F_{n-1\ n})
$$
for the global intertwiner $\Omega$.
It is interesting to note that for such a twist the element
${\cal F}$ carries out an
algebra homomorphism $a\to <{\cal F},id\oo a>$
 from the quantum semi-group
${\cal A}_{R}\sim {\cal H}^*$ to ${\cal H}$, while the transposed
mapping is a homomorphism from
${\cal A}_{\tilde R}$  to $\tilde {\cal H}\sim {\cal H}$. Composition
of these mappings with the representation $\rho$
yields in its turn homomorphisms of the twisted and
non-twisted semi-groups to ${\bf R}$ which are determined on the generators
by the element $F$. The necessary and sufficient conditions for
the existence of
such homomorphisms are just exactly equations (\ref{FC21}--\ref{FC22}).
Thus, there is a tool for verification whether two solutions
to YBE are related via the twist with factorization conditions
(\ref{FC21}) and (\ref{FC22}): among all the invertible elements $F$
intertwining $R$ and $\tilde R$ one should find those defining
homomorphisms from the corresponding quantum groups
into ${\bf R}$.

\section{Conclusion}
The present investigation shows
that transformation $\tilde R= \tau (F^{-1}) R F$ of a solution to YBE
leads to a new solution if there exists a three-tensor $\Omega^3$
relating $R$- and $\tilde R$-symmetric three-tensors. For the global
twist of the bialgebra ${\cal U}$ defined by $R$, one should require
the existence of invertible elements $\Omega^n$ relating $R$- and
$\tilde R$-symmetric n-tensors for every $n$.
This means the equivalence between the corresponding representations of
the braid groups $B_n$. Since dimension $n=3$ proves to be crucial
for the Yang-Baxter equation, an interesting question is
whether two  representations of $B_n$
are locally isomorphic if such an isomorphism takes place for $B_3$.
If so, that could reduce the problem of building twist of bialgebras,
within the matrix formalism, to solving the finite set of relatively
simple equations in the matrix tensor square and cube.

Although twist establishes an equivalence between monoidal categories of
representations of quasitriangular Hopf algebras, the physical content of
related integrable
models can change significantly. So, the jordanian deformation of the
$XXX$-model of spin $\frac{1}{2}$ chain
leads to the non-Hermitian Hamiltonian \cite{KS}. Preservation of its
spectrum under that particular transformation is occasional, rather,
and does not take place in other cases, for example, in transition from
the standard quantum Toda chain to the system related to the
Cremmer-Gervais R-matrix \cite{CKD}. On the other hand, there is a successful
experience of applying twisting technique to obtain simpler expressions
for  correlation functions \cite{MS}, and
the global intertwiner $\Omega$ introduced for the special case
of the $XXZ$-model in \cite{MS}
and studied on a somewhat general  basis in the
present paper should play an essential role in that process.

Another possible application of the present consideration is finding
new solutions to the matrix Yang-Baxter equation including those depending
on the spectral parameter. Particular realization of this line requires
essentially using computer algebra programming because in the simplest
case of two dimensions all the solutions has already been listed in
\cite{Hietarinta}; and that is beyond the scope of our communication
being a separate and elaborate problem.

\newpage
\noindent
{\Large\bf  Acknowledgement}
\vspace{0.3cm}

\noindent
We are grateful to Professor T. Hodges for his valuable comments and
remarks on the subject of this work.

\end{document}